\documentclass{article}
\usepackage{amssymb,latexsym}
\usepackage{amsmath}
\usepackage{graphics}

\newtheorem{theorem}{Theorem}
\newtheorem{corollary}{Corollary}
\newtheorem{lemma}{Lemma}

\numberwithin{equation}{section}
\begin{document}
\title{The distribution of the minimum height among pivotal sites in critical
two-dimensional percolation }

\author{G.J. Morrow\ \ Y. Zhang}

\date{July 21,\ \ 2005.}

\maketitle

\begin{abstract}

\vspace{1pt}Let $L_{n}$ denote the lowest crossing of the $2n\times 2n$ square box $B(n)$ centered at the origin for critical site percolation on $\Bbb{Z}^2$ or critical site percolation on the triangular lattice imbedded in $\Bbb Z^{2}$, and denote by $Q_{n}$ the set of pivotal sites along this crossing. On the event that a pivotal site exists, denote the minimum height that a pivotal site attains above the bottom of $B(n)$ by \[M_{n}:=\min \{m\ge 0: (x,-n+m)\in Q_{n} \mbox{ for some }-n\le x\le n\}\] Else, define $M_{n}=2n$. We prove that 
$P(M_{n} < m)\asymp m/n, \mbox{ uniformly for }1\le m\le n$. This relation extends Theorem 1 of van den Berg and Jarai \cite{BJ03} who handle the corresponding distribution for the lowest crossing in a slightly different context. As a corollary we establish the asymptotic distribution of the minimum height of the set of cut points of a certain chordal $SLE_{6}$ in the unit square of $\Bbb C$.
\end{abstract}

\section{Introduction}\label{S:intro}
     
We shall consider site percolation on the triangular lattice, or site (or bond) percolation on $\Bbb Z^{2}$. Note that site and bond percolation on $\Bbb{Z}^2$ are equivalent for our purposes (see Kesten \cite{K82}, Chp. 3.1).  Our results hold for site percolation on more general planar lattices such that in the critical model, the probability of horizontal or vertical crossings (either open or closed) of an $n$ by $n$ square are bounded below independent of $n$ (see \cite{K82}, Theorem 6.1, or \cite{K87}), but we restrict our notation to handle only the two site percolation models above.  Each vertex of the lattice is open with probability $p$ and closed with probability $1-p$ and the sites are occupied
independently of each other.  We will realize the triangular lattice with vertex set $%
\Bbb Z ^{2}$ as follows. For a given $(x,y)\in \Bbb Z^{2}$, its nearest
neighbors are defined as $(x\pm 1,y)$, $(x,y\pm 1)$, $(x+1,y-1)$, and $%
(x-1,y+1)$. Bonds between neighboring or adjacent sites therefore correspond to
vertical or horizontal displacements of one unit, or diagonal displacements
between two nearest vertices along a line making an angle of $135^{\circ }$
with the positive $x$-axis. Site percolation on $\Bbb Z^{2}$ is defined similarly except the diagonal bonds are excluded. Recall that the triangular lattice may also be viewed with sites as hexagons
in a regular hexagonal tiling of the plane. This point of view is convenient
to describe the fact that when we have critical percolation on the triangular lattice ($p=1/2$) and the
hexagonal mesh tends to zero, the percolation cluster interface has a
conformally invariant scaling limit, namely the Schramm-Loewner evolution
process $SLE_{6}$ (Smirnov \cite{S01b}). It is still unknown whether two-dimensional lattices other than the triangular lattice have this conformally invariant scaling limit. Nevertheless in this paper we will derive our results directly from the percolation structure, so we will not be restricted only to the triangular lattice. In the sequel we will be working exclusively with the critical percolation model.

Define $\left\| \mathbf{x}\right\|
:=\max \{\left| x\right| ,\left| y\right| \}$ for $\mathbf{x}=(x,y)\in 
\Bbb Z ^{2}$. For any real number $r\ge 0$ we denote the square box of vertices $B(r):=\{\mathbf{x}\in \Bbb Z ^{2}:\left\| \mathbf{x}\right\| \leq r\}$. A path is a sequence of distinct vertices connected by nearest neighbor bonds. Thus a path is simple. A horizontal open (closed) crossing of a rectangle $R$ is an open (closed) path in $R$ from the left side of $R$ to the right side of $R$. A vertical crossing is defined similarly. Let $n$ be a positive integer. 
The lowest crossing for any given configuration of vertices for which a horizontal open crossing of $B(n)$ exists is known (see Grimmett \cite{G99}) to be the unique horizontal open crossing $L_{n}$ of $B(n)$ that lies in the region on or beneath any other horizontal open crossing. For each vertex $\mathbf{x}\in L_{n}$ there exist two disjoint open paths from $\mathbf{x}$ to the sides of $B(n)$ and a closed path to the bottom of $B(n)$. This particular orientation of disjoint paths or arms in fact characterizes the vertices of the lowest crossing. Define the set of pioneering sites $F_{n}$ as the union of the lowest crossing with the many complicated orbs and tendrils hanging from it. These latter sets consist of sites $\mathbf{x}$ such that there exists a closed path from $\mathbf{x}$ to the bottom of $B(n)$ and one open path each from $\mathbf{x}$ to the left and right sides of $B(n)$, but such that these open paths are not disjoint. Alternatively, when the lowest crossing exists, $F_{n}$ is the set of open sites connected to this crossing that are discovered through the exploration process that starts at the lower left corner of $B(n)$, and runs until it meets the right side and finally the lower right corner, that determines in particular the interface between the lowest left to right spanning open cluster in $B(n)$ and the closed cluster attaching to its bottom side. In the case of the triangular lattice, this description of $F_{n}$ explains its correspondence to the trace of a certain chordal $SLE_{6}$ (see the description preceding Corollary \ref{C:SLE}).

The pivotal sites are vertices that lie along the lowest crossing and that satisfy the following property: if the vertex is changed from open to closed then there no longer exists a horizontal open crossing. Explicitly, define the event that $\mathbf{x}$ is a pivotal site by
\begin{equation}\label{E:Qxn}
\begin{array}{ll}
{\cal Q}(\mathbf{x},n):=&\mbox{there exists a horizontal open crossing of }B(n)\\ & \mbox{containing the vertex } \mathbf{x},\mbox{ and there exist two }\\ &\mbox{disjoint closed paths in } B(n)\mbox{ started from } \mathbf{x}, \mbox{ one}\\ &\mbox{to the top side and one to the bottom side of }B(n),
\end{array}
\end{equation}
and denote the set of pivotal sites by $Q_{n}:=\{\mathbf{x}: {\cal Q}(\mathbf{x},n) \mbox{ occurs}\}$. We see from the definition that the two closed arms emanating from a pivotal (hence open) site $\mathbf{x}$ force any horizontal open crossing of $B(n)$ to pass through $\mathbf{x}$, and conversely, if $\mathbf{x}$ is open and the pivotal property holds then the definition \eqref{E:Qxn} holds. 

In the triangular lattice case, the lowest crossing corresponds to the upper boundary of an $SLE_{6}$ path in the unit square of $\Bbb{C}$ started from $0$ and aimed toward $1$, while the highest crossing corresponds to the lower boundary of an $SLE_{6}$ path in this square that goes from $1+i$ to $i$.  The pivotal sites that form where the highest and lowest crossings meet correspond to the cut points of $SLE_{6}$ in a sense that we make precise in Corollary \ref{C:SLE} below. Here the outer boundary and cut points of a  process $z(s)$, $s\geq 0$, in $\Bbb{C}$ are defined as follows. Define the hull $K_{t}$ at time $t$ as 
the union of the trace $z [0,t]:=\{z(s):0\leq s\leq t\}$ 
with the bounded components of its complement $\Bbb C \backslash z[0,t].$ 
The frontier or outer boundary of the trace $z[0,t]$ up to time $t$
is defined as the boundary of $K_{t}$.  
By contrast, a pioneer
point is defined as any
point $z(s)$ at some time $s\leq t$ such that $z(s)$ is on the frontier of $z[0,s]$. A point $z(s)$ for some 
$0<s<t$ is called a cut point of $z[0,t]$ if $K_{t}\setminus \{z(s)\}$ is disconnected, so that  
$z(s)$ is in the outer boundary of $z[0,t]$, and $s$ is then called a cut time. Lawler, Schramm, and Werner \cite{LSW01c} have shown that the
frontier, pioneer points, and cut points of a planar Brownian motion almost surely have Hausdorff dimensions respectively $4/3$, $7/4$, and $3/4$. The corresponding dimensions have been obtained for $SLE_{6}$ itself (see Beffara \cite{Beff04} and Lawler, Schramm, and Werner \cite{LSW01a},\cite{LSW01b},\cite{LSW01c}). In the case of $SLE_{6}$ the trace is the same as the set of pioneering points because almost surely the path is continuous and non-simple but non-crossing (Rhode and Schramm \cite{RS03}). Beffara \cite{Beff04} has shown that the dimension of the set of cut times for $SLE_{6}$ is $2/3$. It would be interesting to find a percolation analogue of this result.

In the case of the triangular lattice, Smirnov and Werner \cite{SW01} apply the connection with $SLE_{6}$ to establish estimates of probabilities of $\kappa$-arm paths in both the plane and half-plane. The planar exponents are $(\kappa^2-1)/12$ while  the half-plane exponents appear in Lemma \ref{L:PHE3}. Morrow and Zhang \cite{MZ05} use the planar estimates for $\kappa=2,3,4$ to establish the following discrete analogue of the above Hausdorff fractal dimensions via asymptotic moment relations. If the lattice is triangular, then for all natural numbers $\tau$, as $n \rightarrow \infty$
\begin{equation}\label{E:moments}
\mbox{(i) }E(\left| L_{n}\right|^{\tau}) =n^{4\tau/3+o(1)},\mbox{ \ (ii)}  E(\left| F_{n}\right|^{\tau}) =n^{7\tau/4+o(1)},\mbox{ \ (iii) }
E(\left| Q_{n}\right|^{\tau}) =n^{3\tau/4+o(1)}.
\end{equation}
However due to a $o(1)$ term appearing in the exponents of the planar estimates \cite{SW01}, these estimates are not sufficient to establish Theorem \ref{T:minheightpivotal} below, yet the connection to $SLE_{6}$ is not necessary either. For the proof of Theorem \ref{T:minheightpivotal} we only need to estimate the probability of the occurrence of pivotal sites near the boundary of $B(n)$, so we can rely on exact (no $o(1)$ term), integer-valued half-plane exponents for two- and three-arm paths (Lemma \ref{L:PHE3} below) that can be proved directly in the case of the triangular lattice by change of color arguments (see Lawler, Schramm, and Werner \cite{LSW02b}, Appendix A). In Section \ref{S:prooflemma1} we extend this Lemma \ref{L:PHE3} to more general lattices (without using change of color arguments) by using instead arguments similar to those of Kesten, Sidoravicius, and Zhang \cite{KSZ98}, Lemma 5. These authors establish directly the exact integer exponent $(5^2-1)/12=2$ for a multicolor five-arm exponent in the plane in the case of a two-dimensional lattice. Their arguments rely heavily on the connection method of Kesten \cite{K87} that will also be applied throughout this paper.

Denote the minimum height that a pivotal site attains above the bottom of $B(n)$ when such a site exists by by $M_{n}:=\min \{m\ge 0: (x,-n+m)\in Q_{n} \mbox{ for some }-n\le x\le n\}$, else define $M_{n}=2n$. In this paper we establish the asymptotic distribution of this minimum height as follows.

\begin{theorem}\label{T:minheightpivotal}
$P(M_{n} < m)\asymp m/n, \mbox{ uniformly for }1\le m\le n$.
\end{theorem}
This relation extends Theorem 1 of van den Berg and Jarai \cite{BJ03} who establish a corresponding result for the lowest crossing in a slightly different context. These authors study critical percolation in the closed upper half plane $\overline{\Bbb H}$ and consider the event that the lowest crossing that connects the $x$-axis from $(-\infty,0]$ to $[n,\infty)$ will meet the horizontal line $y=m$ somewhere in the interval $0<x<n$. By using a direct surgery on the configuration space to effect a proof of a certain key near independence relation, 
they find that such a meeting occurs at least once, and therefore the minimum height of the crossing is at most $m$, with asymptotic probability $1/\log (n/m)$.  
To prove our Theorem \ref{T:minheightpivotal}, it turns out that we need only establish asymptotic evaluations of the first and second moments of the variable $X_{n,m}$ defined by \eqref{E:Xnm}. For this we shall consider not only the probability of a two- or three-arm path crossing a semi-annulus, but more generally the probability of a two-arm path crossing an annular sector determined by an angle of aperture $\varphi=\pi/2$. Fortunately such a ``corner" probability can be proved to have a critical exponent strictly larger than the semi-annular case (cf. Kesten and Zhang \cite{KZ87}). 

As a corollary of the proof of Theorem \ref{T:minheightpivotal} combined with the planar estimate of \cite{SW01} for a four-arm path, by following \cite{MZ05}, formula (5.19), we find in the special case of the triangular lattice how the average density of pivotal sites decreases as we go toward the boundary of $B(n)$ . Denote by $Q_{n,m}$ the number of pivotal sites in the horizonatlal strip of height $m$ above the bottom of $B(n)$.
\begin{corollary}\label{C:EQnm}
If the lattice is triangular, then $E(Q_{n,m})=m^{7/4+o(1)}n^{-1}$, as $m\rightarrow \infty$. 
\end{corollary}
Although it is not known how to remove the $o(1)$ term from the exponent in this statement it is possible to make the precise asymptotic statement of Theorem \ref{T:minheightpivotal} not by counting the number of pivotal sites in the bottom of $B(n)$ but by counting only the number of $m$ by $m$ adjacent blocks along the bottom of $B(n)$ that contain a pivotal site and by applying both Lemmas \ref{L:PHE3} and \ref{L:cornerprob} below (see Section \ref{S:momentstrategy}). 

Finally, consider again the critical percolation on the triangular lattice realized with sites as hexagons of diameter $\delta >0$ in a regular hexagonal lattice. We describe the weak convergence of a certain percolation exploration path to $SLE_{6}$ following Camia and Newman \cite{CN05}. Let $D$ be a Jordan domain in $\Bbb C$ (the boundary of $D$ is a continuous simple loop) and let $a$ and $b$ be points on the boundary of $D$. Given $D$, and the boundary points $a$ and $b$, denote by $D^{\delta}$ the largest region of hexagons inside $D$, and define appropriate vertices $\mathbf{x}_{a}$ and $\mathbf{x}_{b}$ of the boundary of $D^{\delta}$ that are nearest to $a$ and $b$. For any configuration $\sigma$ of open and closed sites in $D^{\delta}$, define the exploration path $\gamma_{D,a,b}^{\delta}(\sigma)$ as follows. Introduce closed sites adjacent to and outside the part of the boundary of $D^{\delta}$ that is oriented counterclockwise from $\mathbf{x}_a$ to $\mathbf{x}_b$. Similarly introduce open sites attached to the outside of the remaining part of the boundary of $D^{\delta}$ that is oriented clockwise from $\mathbf{x}_b$ to $\mathbf{x}_a$. The exploration path starts at $\mathbf{x}_{a}$ and proceeds by moving along neigboring sites, keeping closed (hexagonal) sites on its right and open sites on its left until it reaches $\mathbf{x}_{b}$. Define the uniform metric between continuous curves $\gamma_{i}:[0,1]\rightarrow \Bbb{C}$, $i=1,2$, by $d(\gamma_{1},\gamma_{2}):=\inf \sup_{t\in [0,1]}|\gamma_{1}(t)-\gamma_{2}(t)|$ (the infimum is over all reparametrizations).  Smirnov (2001) has shown that the distribution of $\gamma_{D,a,b}^{\delta}$ converges to that of the trace $\gamma_{D,a,b}$ of chordal $SLE_{6}$ inside $D$ from $a$ to $b$ with respect to the uniform metric.  Camia and Newman \cite{CN05} point out that by Billingsley \cite{Bill71}, Corollary 1, the curves $\gamma^{\delta (k)}= \gamma_{D,a,b}^{\delta (k)}$ and $\gamma= \gamma_{D,a,b}$ can be coupled for some sequence $\delta (k) \rightarrow 0$ on some probability space $(\Omega,{\cal B},\Bbb{P})$ such that $d(\gamma^{\delta(k)},\gamma)\rightarrow 0$ for every $\omega \in \Omega$ as $k \rightarrow \infty$. Therefore by the proof shown in Section \ref{S:Cor2} we obtain the following.
\begin{corollary}\label{C:SLE}
Let $D$ be the unit square in $\Bbb C$ with vertices $0$, $1$, $1+i$, and $i$. Let $\gamma=\gamma_{D,1+i,i}$ be the trace of the (upper) chordal $SLE_{6}$ in $D$ from $1+i$ to $i$.  Consider the path $\gamma$ from the last time it leaves the right boundary of $D$ until the first time it meets the left boundary of $D$. Denote by $K$ the lower boundary of this piece of the path $\gamma$. Let ${\cal E}$ be the event that $K$ does not touch the bottom or top of $D$. 
Let ${\cal M}$ denote the minimum height, or smallest vertical distance between a cut point on $\gamma$ and the $x$-axis. Then $P({\cal M}\le s|{\cal E})\asymp s$ as $s\rightarrow 0$.
\end{corollary}

\section{Moment strategy}\label{S:momentstrategy}  
We follow the strategy of van den Berg and Jarai \cite{BJ03}.
For simplicity of the exposition we assume that $n=2^{j_{0}}$ and $m=2^{-i_{0}}n$, for some non-negative integers $j_{0}=j_{0}(n)$ and $i_{0}$. We may easily remove these conditions later. Note that by this choice of parameters, $i_{0}=c\log (n/m)$ and $j_{0}=c\log(n)$. For each $k=1,2,\dots, 2n/m$, we construct the square blocks \[H_{k}:=\{(x,y)\in B(n):-n+(k-1)m<x\le -n+km,-n\le y\le -n+m\}.\] Thus the $H_{k}$ are adjacent $m$ by $m$ blocks of vertices sitting along the bottom of $B(n)$ from $x=-n$ to $x=n$. We count the number of these blocks that contain a pivotal site:
\begin{equation}\label{E:Xnm}
X_{n,m}:=\sum \limits _{k=1}^{2n/m}1_{\{Q_{n}\cap H_{k}\neq \emptyset\}}.
\end{equation}
Notice that the minimum height $M_{n}$ of the set of pivotal sites $Q_{n}$ above the bottom of $B(n)$ is at most $m$ precisely when $X_{n,m}\ge 1$. We then calculate (1) an asymptotic expression for the first moment of $X_{n,m}$, and (2) an upper bound for the second moment of $X_{n,m}$. By the Cauchy-Schwarz inequality this gives the lower bound  
\begin{equation}\label{E:LB}
P(X_{n,m}\ge 1)\ge E(X_{n,m})^2/E(X_{n,m}^2).
\end{equation}
Since in the pivotal case the first and second moments of $X_{n,m}$ turn out both to be asymptotically equal to $m/n$ (up to constants), Markov's inequality will yield an efficient upper bound:
\begin{equation}\label{E:UB}
P(X_{n,m}\ge 1)\le E(X_{n,m}).
\end{equation}
Thus Theorem \ref{T:minheightpivotal} will be proved.

\subsection{Multiple-arm crossings}\label{SS:multiarmcrossings}   Probability estimates for certain horseshoe events will be needed in our computations. We define horseshoe sets of vertices that lie along the bottom of $B(n)$ as follows.
Let $B_{1}\subset B(n)$ be a box such that the bottom edge of $B_{1}$ lies on the bottom edge of $B(n)$, and let $B_{2}\subset B(n)$ be a box containing $B_{1}$ such that the bottom edge of $B_{1}$ is centered in the bottom edge of $B_{2}$. Thus $B_{2}\setminus B_{1}$ is a semi-annular region that we call a horseshoe.  
Let $B_{1}=B_{1}(2^{\rho})\subset B(n)$ be a fixed box of radius $2^{\rho}$ and for each $\nu\ge \rho$, let $B_{2}=B_{2}(2^{\nu})\subset B(n)$ be a box of radius $2^{\nu}$ containing $B_{1}$ such that $(B_{1},B_{2})$ form a horseshoe pair. Denote by $H:=H(\rho,\nu):=B_{2}(2^{\nu})\setminus B_{1}(2^{\rho})$ the corresponding horseshoe. We denote by $\partial_{+}B$ the top and side edges of a box $B$. 
Denote also by ${\cal J}_{\kappa}(\rho,\nu)$ the event that there is a $\kappa$-arm crossing of the horseshoe $H(\rho,\nu)$.  Here if $\kappa=2$ we assume one open path and one closed path oriented counterclockwise from the open path as viewed from the inner horseshoe box, while if $\kappa=3$ we assume the paths are open, closed, open in the counterclockwise order again as viewed from the inner horseshoe box. The paths start from $\partial_{+}B_{1}$ and end on $\partial_{+}B_{2}$ and otherwise remain in $H(\rho,\nu)$. In the following we may even assume if we so choose that the paths start from fixed disjoint intervals along $\partial_{+}B_{1}(2^\rho)$ that have lengths $c2^{\rho}$ for some constant $c>0$. The following Lemma for the probabilities of the events ${\cal J}_{2}(\rho,\nu)$ and ${\cal J}_{3}(\rho,\nu)$ is proved in Section \ref{S:prooflemma1}. There, to avoid the use of ``changing-the-colors" arguments that do not apply to a general two-dimensional lattice, the fence argument of Kesten \cite{K87} is used instead. The integer exponents obtained are the so-called half-plane exponents for the cases $\kappa=2,3$. See Kesten, Sidoravicius, and Zhang \cite{KSZ98} who prove a similar result for a multicolor five-arm exponent in the plane.  

\begin{lemma}\label{L:PHE3} 
For $\kappa = 2,3$ we have $P({\cal J}_{\kappa}(\rho,\nu))\asymp (2^{\rho}/2^{\nu})^{\kappa(\kappa+1)/6}$. 
\end{lemma}
Note that in the case of the triangular lattice, the Lemma \ref{L:PHE3} may be extended to all $\kappa \ge 4$ by using Theorem 3 of \cite{SW01} and Kesten's connection method \cite{K87} at the cost of an additional $o(1)$ term in the exponent as $\rho \rightarrow \infty$ (see \cite{MZ05}, Lemma 5). However it is precisely the asymptotic relation up to constants shown in Lemma \ref{L:PHE3} that is needed for our proof.

\subsection{First moment}\label{SS:firstmoment}
We now calculate the first moment of $X_{n,m}$. Break up the $x$-axis from $0$ to $n$ into dyadic intervals of integers $I_{j}$ with lengths $2^{-j}n$, $j=0,1,2,\dots,i_{0}$. Explicitly we write $I_{j}:=(n-2^{-j}n, n-2^{-j-1}n]$, $j=0,1,\dots,i_{0}-1$ and $I_{i_{0}}=(n-m,n]$. The interval $I_{j}$ is the projection onto the $x$-axis of a certain collection of the adjacent square blocks $H_{k}$. By our choice of $m$ and $n$ this collection consists of $2^{i_{0}-j}$ blocks. Fix for the moment some $n/m\le k\le 2n/m$ and the corresponding interval index $j=j(k)$, where if $(x,y)\in H_{k}$ then $x\in I_{j}$. There is a box $B_{2}=B_{2}(2^\nu)$ with bottom edge along the bottom edge of $B(n)$ and right edge along the right edge of $B(n)$ and with a radius given asymptotically by $2^{\nu}\asymp 2^{-j(k)}n$ such that $H_{k}\subseteq B_{2}(2^\nu) \subset B(n)$ and such that with $B_{1}=H_{k}$ we have a horseshoe pair $(B_{1},B_{2})$ as in the context of Lemma \ref{L:PHE3}. Denote the radius of the inner box $B_{1}$ by $2^{\rho+1}=m$ and denote by ${\cal J}_{3}(\rho,\nu)$ the event of a three-arm crossing of the horseshoe $H(\rho,\nu)=B_{2}\setminus B_{1}$, again as in the context of Lemma \ref{L:PHE3}. Then according to Lemma \ref{L:PHE3} with $\kappa = 3$ we have 
\begin{equation}\label{E:asympthreearm}
P({\cal J}_{3}(\rho,\nu))\asymp (2^{-j}n/m)^{-2}=2^{2j}m^2n^{-2}.
\end{equation}

Now on the event that there exists a pivotal site in $H_{k}$, besides the event ${\cal J}_{3}(\rho,\nu)$, there is a two-arm path from  $\partial_{1}B_{2}(2^\nu)$ to $\partial_{1}B(n)$ where $\partial _{1}B$ denotes the union of the top and left edges of a given box $B$. However since the box $B_{2}(2^\nu)$ lies in the bottom right corner of $B(n)$, the probability of this two-arm path is not governed by the critical exponent of Lemma \ref{L:PHE3} with $\kappa=2$. Indeed the two-arm path is now restricted to lie in a sector with an opening of angle $\varphi =\pi/2$ rather than in a half-plane, so the probability of the two-arm path so restricted has a strictly different critical exponent. To make this statement precise we follow the method of Kesten and Zhang \cite{KZ87}. Define the sector with aperture $\varphi$ inside $B(n)$ by \[S(\varphi,n)=\{(x,y)\in B(n): x=r\cos \theta,y=r\sin \theta\mbox{ for some }r\ge 0,0\le \theta\le \varphi\}.\] Let $\varphi\ge \pi/2$ and let $1\le l<n$. Denote by  $K_{1}:=\{(l,y):3l/8\le y\le 5l/8\}$ an interval of vertices of length $l/4$ on the right edge of $S(\varphi,l)$, and similarly denote by $K_{2}:=\{(x,l):l/4\le x\le 3l/4\}$ a matching interval of vertices along the top edge of $S(\varphi,l)$. These intervals simply sit in the centers of the right and top edges of the square sector $S(\pi/2,l)$. Denote  by $\zeta ( \varphi,l,n)$ the event that there is a closed path in $S(\varphi,n)\setminus S(\varphi,l)$ started from $K_{1}$ and ending on $\partial B(n)$ and also an open path in $S(\varphi,n)\setminus S(\varphi,l)$ started from $K_{2}$ and ending on $\partial B(n)$. 

\begin{lemma}\label{L:cornerprob} There exists a constant $\alpha>0$ such that uniformly for $1\le l\le n/2$, \[P(\zeta(\pi/2,l,n))/P(\zeta(\pi,l,n))\le C(n/l)^{-\alpha}.\] 
\end{lemma}
Lemma \ref{L:cornerprob} is proved for the case of a single arm by \cite{KZ87}, Theorem 2. We show a modification of that proof needed to establish the two-arm case in an Appendix. Note that the method of proof following \cite{KZ87} is not restricted only to comparison of probabilities of two-arm paths between the quarter-plane and half-plane cases, but applies to comparisons between probabilities of two-arm paths for sectors of any different apertures. 

By the connection method of \cite{K87}, Lemma 5 (cf. the fence argument of Section \ref{S:prooflemma1} below), we may argue that the probability that there is a two-arm path in $B(n)\setminus B_{2}$ from $\partial_{1}B_{2}$ to $\partial_{1}B(n)$ is at most $CP(\zeta(\pi/2,2^{\nu+1},n))$. But by Lemma \ref{L:PHE3} with $\kappa=2$ we have $P(\zeta(\pi,2^{\nu+1},n))\asymp 2^{-j}n/n=2^{-j}$. Hence by this precise asymptotic evaluation we obtain by \eqref{E:asympthreearm}, Lemma \ref{L:cornerprob} with $l=2^{\nu+1}$,  and independence that 
\begin{equation}\label{E:ubprobQnHk}
P(Q_{n}\cap H_{k}\neq\emptyset)\le C2^{(1-\alpha)j}m^2n^{-2}.
\end{equation}

We may also obtain a lower bound that will serve an equivalent purpose in our moment strategy to the upper bound \eqref{E:ubprobQnHk} by again employing connection arguments. First by \cite{K87} Lemma 3 (extension of the FKG inequality) one easily constructs (see for example \cite{MZ05}, Lemma 5) by induction a two-arm path in $S(\pi/2,n)\setminus S(\pi/2,l)$ in such a way that  
\begin{equation}\label{E:lbcornerprob}
P(\zeta(\pi/2,l,n))\ge (1/C)(n/l)^{-1-\beta} \mbox{ for some constant }0< \beta < \infty.
\end{equation}
By this method we may even obtain \eqref{E:lbcornerprob} with the added requirement that the two-arm path in $\zeta(\pi/2,l,n)$ meets $\partial_{1}B(n)$ in pre-specified disjoint intervals of length $n/4$. 
Second, we argue that there is a constant $c>0$ such that $P(\exists \mbox{ a pivotal site in }B(r))\ge c$, uniformly in $r\ge 1$. Define certain rectangles that sit inside $B(r)$ and adjacent to the top, left, and right edges of $B(r)$ by  $R_{1}:=[-r/4,r/4]\times [3r/4,r]$, $R_{2}:=[-r,-3r/4]\times [-r,r/4]$, and $R_{3}:=[3r/4,r]\times [-r,r/4]$. Let $R$ be a rectangle with sides parallel to the coordinate axes and sharing one side with the boundary of a box $B$. We say that a path $h$-tunnels through $R$ on its way to $\partial B$ if the intersection of the path with the smallest infinite vertical strip containing $R$ remains in $R$. Thus the path may weave in and out of $R$ but not through the top or bottom sides of $R$, and comes finally to $\partial B$. Likewise we say that a path $v$-tunnels through $R$ on its way to $\partial B$ if the roles of horizontal and vertical are interchanged in the preceding definition. This definition is consistent with the requirements of \cite{K87}, Lemma 4. Define also the rectangle $R_{h}(r):=[-r,r]\times[-r,r/4]$ that spans $B(r)$ from left to right and sits on the bottom of $B(r)$, and the rectangle $R_{v}(r):=[-r/4,r/4]\times [-r,-r]$ that spans $B(r)$ from bottom to top. Now, first since we have critical percolation, by the RSW theory there is a constant $c_{1}>0$ such that $P(\exists \mbox{ a horizontal crossing of }R_{h}(r))\ge c_{1}$. We may decompose the event that there exists a horizontal crossing of $R_{h}(r)$ by specifying a lowest horizontal crossing ${\cal L}$ of $R_{h}(r)$. Conditional on the event that the lowest crossing is ${\cal L}=\gamma_{0}$ for some fixed path $\gamma_{0}$, we have independence between the configuration of sites above $\gamma_{0}$ and those below it. There is a constant $c_{2}>0$ such that with probability at least $c_{2}$ on the configuration of sites above $\gamma_{0}$ there exists a closed path in the rectangle $R_{v}(r)$ that meets $\gamma_{0}$ on the bottom and also meets the top of $B(n)$. Therefore by summing over $\gamma_{0}$ we have that $P(\exists \mbox{ a pivotal site in }B(r))\ge c_{1}c_{2}$. Moreover we can construct vertical open crossings of $R_{2}$ and $R_{3}$ and a horizontal closed crossing of $R_{1}$ such that by independence and FKG as in Lemma 3 of \cite{K87} we have in fact that there exists a constant $c_{3}>0$ such that with probability at least $c_{3}$, in addition to the existence of the pivotal site, these additional crossings exist, and the lowest crossing ${\cal L}$ of $R_{h}(r)$ $h$-tunnels through each of $R_{2}$ and $R_{3}$ and the closed path to the top of $B(r)$ $v$-tunnels through $R_{1}$. We shall apply this construction now with the $m$ by $m$ box $H_{k}$ in place of $B(r)$ and all the rectangles above appropriately scaled and translated to fit inside $H_{k}$ similar as they fit inside $B(r)$.

Finally, let $R_{2,1}(\nu)$, $R_{2,2}(\nu)$, and  $R_{2,3}(\nu)$ be squares with sides parallel to the coordinate axes and with sides of length $2^{\nu -1}$ such that these squares are bisected respectively by the center intervals (of length $2^{\nu-1}$) of the left, top, and right edges of $B_{2}(2^{\nu})$. Note that the right half of the square $R_{2,3}(\nu)$ sits outside the box $B(n)$ by our construction of the corner box $B_{2}(2^{\nu})$. Similarly, construct also squares $R_{1,i}(\rho)$, $i=1,2,3$, that have side lengths $2^{\rho-1}$ such that these squares are bisected respectively by the center intervals of the left, top, and right sides of $B_{1}(2^{\rho})$. The halves of these latter squares that sit inside the inner horseshoe box $B_{1}(2^{\rho})=H_{k}$ coincide respectively with the rectangle $R_{1}$ and the top two-fifths the rectangles $R_{2}$ and $R_{3}$ of the previous paragraph as applied there with $H_{k}$ in place of $B(r)$. By application of Lemma \ref{L:PHE3} to the horseshoe $H(\rho,\nu)$, we may assume that there exists a three-arm path that is open, closed, open in clockwise order where the three paths issue respectively from the center intervals of the left, top, and right edges of $B_{1}(2^{\rho})$. Now by the method of \cite{K87} Lemmas 2 and 4 (see the fence argument of Section \ref{S:prooflemma1}), we may assume that the event ${\cal J}_{3}(\rho,\nu)$ appearing in the conclusion of Lemma \ref{L:PHE3} satisfies the following additional requirements: The closed path is extended such that it $v$-tunnels through $R_{2,1}(\nu)$ on its way to the top side of $R_{2,1}(2^{\nu})$. Further the open path that is counterclockwise from the closed path in $H(\rho,\nu)$ is extended such that it $h$-tunnels through $R_{2,2}(\nu)$ on its way to the left side of $R_{2,2}(\nu)$. Moreover there are horizontal closed crossings of each of the top and bottom halves of the square $R_{2,1}(\nu)$, there is a vertical closed crossing of $R_{2,1}(\nu)$, there is a horizontal closed crossing of the top half of $R_{1,1}(\rho)$, and there is a vertical closed crossing of $R_{1,1}(\rho)$. Finally there are vertical open crossings of each of the left and right halves of the square $R_{2,2}(\nu)$, there is a horizontal open crossing of $R_{2,2}(\nu)$, there is a vertical open crossing of the each of the outer halves of $R_{1,2}(\rho)$ and $R_{1,3}(\rho)$, and there are horizontal open crossings of each of the squares $R_{1,2}(\rho)$ and $R_{1,3}(\rho)$. Hence by the Lemma 3 of \cite{K87} (extended FKG) and the construction of the four-arm path in the previous paragraph, we obtain by the above connection constructions, Lemma \ref{L:PHE3} with $\kappa =3$, and an application of  \eqref{E:lbcornerprob} with $l=2^{\nu}$, that 
\begin{equation}\label{E:lbprobQnHk}
P(Q_{n}\cap H_{k}\neq\emptyset)\ge (1/C)c_{3}(2^{-j}n/m)^{-2}2^{(-1-\beta )j}\ge C2^{(1-\beta)j}m^2n^{-2}.
\end{equation}

We are now ready to sum up our probability bounds to calculate the first moment of $X_{n,m}$. Since there are $2^{i_{0}-j}$ blocks $H_{k}$ in $I_{j}$, and since the interval index $j\le i_{0}$ we have by \eqref{E:ubprobQnHk} that
\begin{equation}\label{E:firstmoment}
E(X_{n,m})\le C\sum \limits_{j=0}^{i_{0}}2^{i_{0}-j}2^{(1-\alpha)j}m^2n^{-2} \le C2^{i_{0}}m^2n^{-2}=Cm/n.
\end{equation}
In the same way but now by \eqref{E:lbprobQnHk}, we have  $E(X_{n,m})\ge Cm/n$. Therefore
\begin{equation}\label{E:asympfirstmoment}
E(X_{n,m})\asymp m/n.
\end{equation}

\subsection{Second moment}\label{SS:secondmoment}
We next calculate the second moment of $X_{n,m}$. We must calculate $P(Q_{n}\cap H_{k_{1}}\neq\emptyset,Q_{n}\cap H_{k_{2}}\neq\emptyset)$ for $k_{2} > k_{1}$. Let $k_{1}\in I_{j_{1}}$ and $k_{2}\in I_{j_{2}}$ for some $0\le j_{1}\le j_{2}\le i_{0}$. There are two cases to consider, (i) $0\le j_{2}-j_{1}\le 4$, and (ii) $j_{2}-j_{1}\ge 5$. Consider case (i) first. We denote the $x$-coordinates of the lower left vertices of the boxes $H_{k_{1}}$ and $H_{k_{2}}$ respectively by $x_{1}$ and $x_{2}$. We then construct estimates based on the size of $x_{2}-x_{1}$. We must have $m\le x_{2}-x_{1}< 2^{j_{1}+1}n$. Fix $x_{1}\in I_{j_{1}}$ and say $x_{2}\in J_{i}$ if $2^{-i}n\le x_{2}-x_{1}< 2^{i+1}n$ for some $i=j_{1},j_{1}+1,\dots,i_{0}$. If $x_{2}\in J_{i}$ then we construct two disjoint boxes $B_{2,k_{1}}(2^{\nu_{1}})$ and $B_{2,k_{2}}(2^{\nu_{1}})$ each with a radius $2^{\nu_{1}}\asymp 2^{-i}n$ such that $(H_{k_{1}},B_{2,k_{1}}(2^{\nu_{1}}))$ and $(H_{k_{2}},B_{2,k_{2}}(2^{\nu_{1}}))$ each form a horseshoe pair. Further by (i) there is a box $B_{1}(2^{\rho})$ containing both these horseshoes and sitting inside $B(n)$ with radius $2^{\rho}\asymp 2^{-i}n$, and a corresponding box $B_{2}(2^{\nu})$ with radius $2^{\nu}\asymp 2^{-j_{1}}n$ that sits in the lower right corner of $B(n)$  such that $(B_{1}(2^{\rho}),B_{2}(2^{\nu}))$ also form a horseshoe pair. On the event that there is a pivotal site in each of $H_{k_{1}}$ and $H_{k_{2}}$, we have that there exist three-arm crossings of each of the three disjoint horseshoes thus constructed. Further there is a two arm-path in $B(n)\setminus B_{2}(2^{\nu})$ from $\partial_{1}B_{2}(2^{\nu})$ to $\partial_{1}B(n)$. Now by the method of \cite{K87} Lemma 5, we may estimate the probability of this two-arm path from above by $CP(\zeta(\pi/2,2^{\nu},n))$. Thus by the same method applied earlier to estimate $P(Q_{n}\cap H_{k}\neq \emptyset)$ from above but this time with Lemma \ref{L:PHE3} applied three times with $\kappa =3$, we have in case (i) that
\begin{equation}\label{E:ubprobQnHk1Hk2i}
P(Q_{n}\cap H_{k_{1}}\neq\emptyset,Q_{n}\cap H_{k_{2}}\neq\emptyset)\le C(2^{-i}n/m)^{-4}2^{-2i+2j_{1}}2^{(-1-\alpha)j_{1}}=C(m/n)^42^{2i+(1-\alpha)j_{1}}.
\end{equation}
Denote now $E(X_{n,m}^2)=\Sigma_{1}+\Sigma_{2}$, where $\Sigma_{1}$ and $\Sigma_{2}$ are the respective sums of the joint probabilities $P(Q_{n}\cap H_{k_{1}}\neq\emptyset,Q_{n}\cap H_{k_{2}}\neq\emptyset)$ over the pairs of indices $(k_{1},k_{2})$ under cases (i) and (ii), where we include now in addition the case $k_{1}=k_{2}$ (corresponding to the sum $E(X_{n,m})$) in the sum $\Sigma_{1}$. Since there are on the order of $2^{i_{0}-i}$ indices $k_{2}$ with $x_{2}\in J_{i}$ for $x_{1}$ fixed in $I_{j_{1}}$, we have by \eqref{E:ubprobQnHk1Hk2i} that
\begin{equation}\label{E:ubSigma1}
\Sigma_{1}\le C(m/n)^4 \sum \limits _{j_{1}=0}^{i_{0}}\sum \limits _{i=j_{1}}^{i_{0}}2^{i_{0}-j_{1}}2^{i_{0}-i}2^{2i+(1-\alpha)j_{1}} \le Cm/n.
\end{equation}
By the lower bound for $E(X_{m,n})$ we have trivially that $\Sigma_{1}\ge (1/C)m/n$. To obtain an asymptotic evaluation of $E(X_{n,m}^2)$ it therefore remains only to show that $\Sigma_{2}\le Cm/n$. To do this we now construct in case (ii) disjoint boxes $B_{2,k_{1}}(2^{\nu_{1}})$ and $B_{2,k_{2}}(2^{\nu_{2}})$ with respective radii $2^{\nu_{1}}\asymp2^{-j_{1}}n$ and $2^{\nu_{2}}\asymp 2^{-j_{2}}n$ , such that again we have horseshoe pairs $(H_{k_{1}},B_{2,k_{1}}(2^{\nu_{1}}))$ and $(H_{k_{2}},B_{2,k_{2}}(2^{\nu_{2}}))$ but now with different outer radii in $B(n)$. This time by (ii) there exists an extra corner box $B'(2^{\nu'})$ with radius $2^{\nu'}\asymp 2^{-j_{1}}n$ that contains the smaller horseshoe associated with the block $H_{k_{2}}$ but that is still disjoint from the larger horseshoe associated with the block $H_{k_{1}}$. Finally there is a second corner box $B''(2^{\nu''})$ with radius $2^{\nu''}\asymp 2^{-j_{1}}n$ that contains both the horseshoes already constructed. Therefore by the same method of estimation as shown for the case (i) we have in case (ii) that
\begin{equation}\label{E:ubprobQnHk1Hk2ii}
\begin{array}{ll}
P(Q_{n}\cap H_{k_{1}}\neq\emptyset,Q_{n}\cap H_{k_{2}}\neq\emptyset)&\le C(m/n)^42^{2(j_{1}+j_{2})}2^{(-1-\alpha)(j_{2}-j_{1})}2^{(-1-\alpha)j_{1}}\\ &\le C(m/n)^42^{2j_{1}+(1-\alpha)j_{2}}.
\end{array}
\end{equation}
Therefore by \eqref{E:ubprobQnHk1Hk2ii},  
\begin{equation}\label{E:ubSigma2}
\Sigma_{2}\le C(m/n)^4 \sum \limits _{j_{1}=0}^{i_{0}}\sum \limits _{j_{2}=j_{1}}^{i_{0}}2^{i_{0}-j_{1}}2^{i_{0}-j_{2}}2^{2j_{1}+(1-\alpha)j_{2}} \le C(m/n)^2\sum \limits _{j_{1}=0}^{i_{0}} 2^{(1-\alpha)j_{1}}\le C(m/n)^{1+\alpha}.
\end{equation}
Hence we have that $\Sigma_{2}$ is in fact of smaller order than $\Sigma_{1}$. Therefore by \eqref{E:ubSigma1} and \eqref{E:ubprobQnHk1Hk2ii} we have established that 
\begin{equation}\label{E:asympsecondmoment}
E(X_{n,m}^2)\asymp m/n.
\end{equation}

The proof of Theorem \ref{T:minheightpivotal} follows now from \eqref{E:asympfirstmoment} and \eqref{E:asympsecondmoment} by using the Cauchy-Schwarz inequality \eqref{E:LB} and Markov's inequality \eqref{E:UB}.$\Box$
\section{Proof of Corollary \ref{C:SLE}}\label{S:Cor2} 
It turns out that the proof we gave of Theorem \ref{T:minheightpivotal} in Section \ref{S:momentstrategy} extends also to the minimum height $M_{hc}$ of the sites on the highest open horizontal crossing of the box $B(n)$. For this reason the main difficulty of Corollary \ref{C:SLE} is to show a lower bound for the distribution of ${\cal M}$. In fact on the event ${\cal E}$ the highest open horizontal crossing will exist with high probability. But the highest horizontal open crossing lies on or below the exploration path. Thus by the approximation of the $SLE_{6}$ path by the exploration path, the upper bound for the conditional distribution of ${\cal M}$ given ${\cal E}$ will follow easily by using the distribution of the minimum height of the highest horizontal open crossing. 
We assume the coupling of exploration paths $\gamma^{\delta(k)}=\gamma^{\delta(k)}_{D,1+i,i}$and $SLE_{6}$ path $\gamma=\gamma_{D,1+i,i}$ on some probability space $(\Omega,{\cal B},\Bbb{P})$ mentioned before. However for convenience we will continue to designate probability by $P$ in place of $\Bbb{P}$ in what follows. Let $z_{\delta(k)}\in \Bbb {C}$ denote the site where the minimal height of a pivotal site is attained within the hexagonal lattice in $D^{\delta(k)}$ whenever such a pivotal site exists, else we let $z_{\delta(k)}=1+i$.  
In the same way that we defined $M$, but now after scaling, define $M^{\delta(k)}$ as the height of $z_{\delta(k)}$.  Denote also the minimum height of the highest horizontal open crossing of the hexagonal lattice configuration in $D^{\delta(k)}$ by $M^{\delta(k)}_{hc}$ when this crossing exists, else we set $M^{\delta(k)}_{hc}=1$.  Let $\epsilon >0$. Denote by ${\cal K}_{\epsilon}$ the event that the maximum height of the upper $ SLE_{6}$ lower boundary $K$ is greater than $1-\epsilon$ or the minumum height is less than $\epsilon$.  By compactness and the contnuity of $\gamma$ \cite{RS03}, we know that $P( {\cal K}_{\epsilon}|{\cal E})\rightarrow 0$ as $\epsilon \rightarrow 0$. 
Therefore by the coupling of $\gamma^{\delta(k)}$ and $\gamma$ we have that the conditional probability of the existence of a horizontal open crossing can be made as close to 1 as we like by first taking $\epsilon$ small and then $k$ large. Note that by these considerations and the fact that there is a positive constant $c$ such that the probability that there exists a horizontal open crossing is bounded below by $c$ uniformly in $k$ imply in particular that $P({\cal E})>0$. Furthermore, if now $1>s>\epsilon>0$ then by our previous observations in this paragraph, we have that
\begin{equation}\label{E:upperboundCor2b}
P({\cal M}\le s, {\cal E})\le P(d(\gamma^{\delta(k)},\gamma)>\epsilon)+P( M^{\delta(k)}_{hc}<s+\epsilon)+P( {\cal K}_{\epsilon},{\cal E}).
\end{equation}
By the analogue of Theorem \ref{T:minheightpivotal} for the highest horizontal open crossing this yields the upper bound for Corollary \ref{C:SLE}.
 
We now turn to the lower bound. In order to force the highest horizontal open crossing in the hexagonal lattice of level $\delta(k)$ to stay away from the top of the domain $D^{\delta(k)}$, we define ${\cal S}_{k}$ as the event that there is a closed horizontal crossing of the horizontal strip $3/4\ge\Im {z}\ge 1/2$.  Let $1/2>s>0$ and choose now $0<\epsilon<s^2/2$.  We further define the following events.
\begin{equation}\label{E:events}
\begin{array}{ll}
&{\cal A}_{k}:= \{ M^{\delta(k)}\le s\}\cap {\cal S}_{k},\mbox{ \ }{\cal D}_{k}:=\{d(\gamma^{\delta(k)},\gamma)\le s^2/2\},\mbox{ \ }{\cal H}_{k}:=\{ M^{\delta(k)}_{hc}\ge s^2\}\\ &{\cal N}_{k}:=\{\exists \mbox{ cut point }z\mbox{ with }|z-z_{\delta(k)}|\le s\}.
\end{array}
\end{equation}
We note that $P( {\cal D}^{c}_{k})=o(1)$ as $k \rightarrow \infty$ by the coupling referred to above. 
On the event ${\cal A}_{k}\cap {\cal H}_{k}$ we have that the highest open horizontal crossing exists and further, the lowest sites on this crossing are at least a vertical distance $s^2$ away from the bottom of $D^{\delta(k)}$, and the exploration path stays below the level $\Im z =3/4$. Since the highest open crossing lies on or below the exploration path in $D^{\delta(k)}$, we have on ${\cal D}_{k}\cap {\cal H}_{k}$ that the lower boundary $K$ must remain a vertical distance $s^{2}/2$ above the bottom of $D^{\delta(k)}$. Therefore, for large $k$ we have that \[ {\cal A}_{k}\cap {\cal D}_{k}\cap {\cal H}_{k}\subset {\cal K}_{s^{2}/4}\subset{\cal E}.\]   Consequently we may estimate that 
\begin{equation}\label{E:lowerboundcalM}
\begin{array}{ll}
P({\cal M}\le 2s, {\cal E})&\ge  P({\cal A}_{k}\cap {\cal D}_{k}\cap {\cal H}_{k}\cap {\cal N}_{k})\\ &\ge P({\cal A}_{k})-P({\cal H}^{c}_{k})-P({\cal D}^{c}_{k})-P({\cal A}_{k}\cap  {\cal H}_{k}\cap {\cal N}^{c}_{k}).
\end{array}
\end{equation} 
By applying the proof of \eqref{E:lbprobQnHk} to the bottom half of $B(n)$ and by then applying the extended FKG inequality (Lemma 3 in \cite{K87}) as before, we may easily show that \eqref{E:lbprobQnHk} continues to hold when we replace $P(Q_{n}\cap H_{k}\neq\emptyset)$ by $P(Q_{n}\cap H_{k}\cap S_{n}\neq\emptyset)$, where $S_{n}$ is the event that there is a closed horizontal crossing of $B(n)$ in the horizontal strip $n/2\le y\le 3n/4$. 
Hence the distribution of $M^{\delta(k)}$ conditional on ${\cal S}_{k}$ is also given as in Theorem \ref{T:minheightpivotal}. Thus since $P(M^{\delta(k)}\le s,{\cal S}_k)\ge C_{1}s$ and $P(M^{\delta(k)}_{hc}\le s^2)\le C_{2}s^2$ we can thereby estimate the first two terms on the right side of \eqref{E:lowerboundcalM}, and we will obtain the lower bound we desire as long as the fourth term can be shown to be of order $o(1)$ as $k \rightarrow \infty$.
We turn therefore to an upper bound for $P({\cal A}_{k}\cap{\cal H}_{k}\cap {\cal N}^{c}_{k})$. Let $0<\epsilon'<\epsilon/10$. Define further events
\begin{equation}\label{E:Gepsilon}
\begin{array}{ll}
{\cal G}_{\epsilon,\epsilon',k}:=&\exists t \mbox{ with }|\gamma(t)-z_{\delta(k)}|\le \epsilon', \mbox{ \ } d(\gamma^{\delta(k)},\gamma)<\epsilon', \mbox{ and }\\ &\exists t_{1}<t<t_{2} \mbox{ with }\gamma(t_{1})=\gamma(t_{2})\mbox{ and  }|\gamma(t_{1})-\gamma(t)|>\epsilon, \mbox{ and }\\ &s^2\le \Im{z_{\delta(k)}}\le s,\mbox{ \ }\epsilon\le \Re{z_{\delta(k)}}\le 1-\epsilon  
\end{array}
\end{equation}
and
\begin{equation}\label{E:Oepsilon}
\begin{array}{ll}
{\cal O}_{\epsilon}:= &\exists t,\exists t_{1}<t<t_{2} \mbox{ with }\gamma(t_{1})=\gamma(t_{2}),\mbox{ and  }\\ &|\gamma(t_{1})-\gamma(t)|\le\epsilon,\gamma(t_{1})\neq \gamma(t_{2})
\end{array}
\end{equation}
and 
\begin{equation}\label{E:Xepsilon}
{\cal X}_{\epsilon,k}:=\{M^{\delta(k)}\le s, \Re z_{\delta(k)}\ge 1-\epsilon\mbox{ or }\Re z_{\delta(k)}\le \epsilon\}.
\end{equation} 
We have 
\begin{equation}\label{E:estPAkNkc}
P({\cal A}_{k}\cap{\cal H}_{k}\cap {\cal N}^{c}_{k})\le P({\cal G}_{\epsilon,\epsilon',k})+P(d(\gamma^{\delta(k)},\gamma)\ge\epsilon')+P({\cal O}_{\epsilon})+P({\cal X}_{\epsilon,k}).
\end{equation}
We comment that there is zero probability that there would be a triple point for the $SLE_{6}$ path, namely \[P(\exists t,\exists t_{1}<t<t_{2} \mbox{ with }\gamma(t_{1})=\gamma(t_{2})=\gamma(t))=0.\] This will actually come out in the proof method below that relies on a six-arm argument. Indeed the triple point for $SLE_{6}$ would imply nine-arms for the exploration path. By the almost sure continuity of $\gamma$ \cite{RS03}, we have that $P({\cal O}_{\epsilon})\rightarrow 0$ as $\epsilon \rightarrow 0$. Further, for fixed $\epsilon'$, we have that $P(d(\gamma^{\delta(k)},\gamma)\ge \epsilon')\rightarrow 0$ as $k \rightarrow \infty$. Finally by defining $X_{\epsilon,k}$ as the number of pivotal sites in the region $\{z:\Im z\le s,\Re z\ge 1-\epsilon\mbox{ or }\Re z\le \epsilon\}$ of $D^{\delta(k)}$, we have $P({\cal X}_{\epsilon,k})=P(X_{\epsilon,k}\ge 1)$. But by the same calculations as shown in Section \ref{S:momentstrategy} (see \eqref{E:firstmoment}) we have that $E(X_{\epsilon,k})\le C(\epsilon)s$ where $C(\epsilon)\rightarrow 0$ as $\epsilon \rightarrow 0$. Thus it will suffice to estimate  $P({\cal G}_{\epsilon,\epsilon',k})$.

Now on the event ${\cal G}_{\epsilon,\epsilon',k}$ we have that the exploration path is within $\epsilon'$ of the $SLE_{6}$ path $\gamma$ in the uniform metric. Indeed that is why there must be some point $\gamma(t)$ within $\epsilon'$ of $z_{\delta(k)}$ since this pivotal site itself must belong to the exploration path. However also on ${\cal G}_{\epsilon,\epsilon',k}$ there is a loop for $\gamma$ that is closed at $\gamma(t_{1})$ with $|\gamma(t_{1})-\gamma(t)|>\epsilon$. Now follow the exploration path from the pivotal site either forward or backward in time corresponding to whether the point of the exploration path approximating $\gamma(t_{1})$ in the uniform metric is either forward in time or backward in time from the time associated with this pivotal site, in a parametrization that makes the uniform metric between the exploration path and the $SLE_{6}$ path at most $\epsilon'$. Assume without loss of generality that this argument leads to following the exploration path forward in time. We follow the exploration path from the pivotal site until it reaches a point $z_{1}$ within $\epsilon'$ of $\gamma(t_{1})$. Then the exploration path comes back to a point within $2\epsilon'$ of the pivotal site as it nears $\gamma(t)$. Finally, proceeding still forward in time, the exploration path comes again within $\epsilon'$ of $\gamma(t_{2})=\gamma(t_{1})$. Therefore there must exist a six-arm path in the percolation configuraion crossing the annulus that is the region outside the  ball $b(z_{\delta(k)},2\epsilon')\subset \Bbb{C}$ (with center the pivotal site and radius $2\epsilon'$) and inside the concentric ball $b(z_{\delta(k)},\epsilon-\epsilon')$. We argue this as follows. Let ${\cal C}_{+}^{*}$ be the cluster of closed sites that issues from the pivotal site $z_{\delta(k)}$ and reaches the top of $D^{\delta(k)}$, and similarly let ${\cal C}_{-}^{*}$ be the cluster of closed sites that connects the pivotal site to the bottom of $D^{\delta(k)}$. Let $r_{1}^{*}$ be the clockwise-most path of closed sites that belongs to the cluster ${\cal C}_{+}^{*}$ such that $r_{1}^{*}$ issues from $z_{\delta(k)}$ and reaches the boundary of $D^{\delta(k)}$. Similarly let $r_{3}^{*}$ be the clockwise-most path of closed sites that belongs to ${\cal C}_{-}^{*}$ that issues from $z_{\delta(k)}$ and reaches the boundary of $D^{\delta(k)}$. Assume without loss of generality but for ease of description that in fact the point $z_{1}$ on the exploration path that we are following as we near $\gamma(t_{1})$ is on the right side of the ball $b(z_{\delta(k)},\epsilon-\epsilon')$ as it is divided by $r_{1}^{*}\cup r_{3}^{*}$. First two long arms (they go all the way to the boundary of $D^{\delta(k)}$) are contributed as follows: one open arm $r_{2}$ along the highest open horizontal crossing to the left side of $r_{1}^{*}\cup r_{3}^{*}$, and one closed arm $r_{3}^{*}$. Another (shorter) open arm $r_{4}$ exists to the right of $r_{1}^{*}\cup r_{3}^{*}$ as the initial part of the exploration path just described from $z_{\delta(k)}$ until we reach $z_{1}$.  Also since this arm can not be traversed backward, it implies an open arm $r_{6}$ being the continuation of the right side of the exploration as it goes on from $z_{1}$ until it tracks forward in time to $\partial b(z_{\delta(k)},2\epsilon')$. Now on the left of the exploration path that has $r_{6}$ on its right there must be a closed path $r_{5}^{*}$. Therefore $r_{6}$ and $r_{5}^{*}$ both lead from $\partial b(z_{\delta(k)},\epsilon-\epsilon')$ to $\partial b(z_{\delta(k)},2\epsilon')$. Now because as the exploration path goes to the right of the pivotal site it must have a path of closed sites on its left, and since this closed path must also connect to the cluster ${\cal C}_{+}^{*}$, we have in fact that this path of closed sites lies to the left of $r_{4}\cup r_{6}$. But since $r_{1}^{*}$ must leave the ball $b(z_{\delta(k)},\epsilon)$, we argue that if  $r_{5}^{*}$ intersects $r_{1}^{*}$ in the region outside $b(z_{\delta(k)},2\epsilon')$, then it must simply be a section of $r_{1}^{*}$ as it continues forward in time past the point $z_{1}$ since otherwise the full path $r_{1}^{*}$ will cease to be a simple path before it leaves $b(z_{\delta(k)},\epsilon)$.   Therefore we have two more closed paths: $r_{5}^{*}$ and either (a) the initial part of $r_{1}^{*}$ until it reaches $z_{1}$ if $r_{5}^{*}$ meets $r_{1}^{*}$ as above, or (b) the full $r_{1}^{*}$, and two open paths: $r_{4}$ and $r_{6}$. Note that in case (a) we actually have seven arms crossing the annulus.

To finish the proof we follow a standard six-arms argument as described in the proof of Lemma 5.1 in Camia and Newman \cite{CN05}. We know from an a priori estimate (cf. \cite{KSZ98}) that the probability of the six-arm path is at most $C(\epsilon'/\epsilon)^{2+\alpha}$ for some constant $\alpha >0$. Since the six-arm event will imply a six-arm crossing of at least one of a grid of approximately $(\epsilon')^{-2}$ annuli with inner radii $4\epsilon'$ and outer radii $(\epsilon-\epsilon')/4$, we obtain $P({\cal G}_{\epsilon,\epsilon',k})\le C(\epsilon')^{-2}(\epsilon'/\epsilon)^{2+\alpha}$. Hence for given $\epsilon>0$ we choose $\epsilon'$ so small that $P({\cal G}_{\epsilon,\epsilon',k})<\epsilon$ and then $k=k(\epsilon)$ so large that $P(d(\gamma^{\delta(k)},\gamma)>\epsilon')<\epsilon$. Therefore by \eqref{E:lowerboundcalM} and \eqref{E:estPAkNkc} we have 
completed the proof of the lower bound for Corollary \ref{C:SLE}. $\Box$

\section{Proof of Lemma \ref{L:PHE3}.}\label{S:prooflemma1}
To prove Lemma \ref{L:PHE3} we work first with the case $\kappa=2$. The proof of the case $\kappa=3$ will follow similar lines. Assume first that $\rho=0$. We assume that $\nu$ is an integer since $P({\cal J}_{2}(0,\nu))$ is decreasing in $\nu$. Let $H$ be a box with sides parallel to the coordinate axes, and let $\Gamma_{H}$ denote the portion of the left edge of this box that is centered within and half as long as this edge. For ease of notation we write $n=2^{\nu}$. Denote by$H(n):=[0,2n]\times [-n,n]=(n,0)+B(n)$ a square box of radius $n$ that sits in the right half plane with left edge along the $y$-axis. We have $\Gamma_{H(n)}=\{(0,y):-n/2\le y\le n/2\}$.  For each $\mathbf{y}\in \Gamma_{H}$, we define the event \[{\cal L}(\mathbf{y},H):= \exists \mbox{ a lowest horizontal open crossing of }H \mbox{ started from }\mathbf{y}\mbox{ on the left edge of }H.\] We want to show that 
\begin{equation}\label{E:sizekappa2}
P({\cal L}(\mathbf{y},H(n))) \asymp 1/n,\mbox{ uniformly for }\mathbf{y}\in \Gamma_{H(n)}. 
\end{equation}
The basic idea is to show that there exists a constant factor $C$ such that for any $\mathbf{y}\in \Gamma_{H(n)}$ and $\widetilde{\mathbf{y}}\in \Gamma_{H(16n)}$, we have
\begin{equation}\label{E:comparePLyB}
P({\cal L}(\mathbf{y},H(n)))\le CP({\cal L}(\widetilde{\mathbf{y}},H(16n))).
\end{equation}
To prove \eqref{E:comparePLyB}, we apply the method of proof of Lemma 4 in \cite{K87}. Denote by $S_{1}(N):=[7N/8,9N/8]\times [0,N/4]$ and $S_{2}(N):=[7N/4,2N]\times [-N/8,N/8]$ squares of radii $N/8$ that sit inside and adjacent to the bottom and right sides of $H(N)$, respectively. Let $\widetilde{{\cal L}}(\mathbf{y},H(N))$
be the event that ${\cal L}(\mathbf{y},H(N)))$ occurs and further there exists a closed path in $H(N)$ to the bottom of $H(N)$ started from $\mathbf{y}$ that will $v$-tunnel through $S_{1}(N)$ and there exists an open path in $H(N)$ to the right side of $H(N)$ started from $\mathbf{y}$ will $h$-tunnel through $S_{2}(N)$, and in addition there is a horizontal closed crossing of $S_{1}(N)$ and a vertical open crossing of $S_{2}(N)$. 

By the method of proof of Lemmas 2 and 4 in \cite{K87} we show that there exists a constant $C_{1}$ such that
\begin{equation}\label{E:comparePLyBPtildeLyB}
P({\cal L}(\mathbf{y},H(n)))\le C_{1}P(\widetilde{{\cal L}}(\mathbf{y},H(4n))), \mbox{ uniformly for }\mathbf{y}\in \Gamma_{H(n)}.
\end{equation} 
Indeed we modify slightly Kesten's approach \cite{K87} to handle the case that the point $\mathbf{y}$ where a two-arm path issues (or, in the case $\kappa=3$, a three-arm path issues) is nearly centered on the left side of $H(n)$ rather than near the center of $H(n)$. Denote by $H(\mathbf{y},n):=\mathbf{y}+H(n)$, the translation of $H(n)$ by $\mathbf{y}$. Also denote by $\partial ^{+}H$ the union of the top, right, and bottom edges of a box $H$. We define an $(\eta,\nu)$-fence for small $\eta>0$ similar as in \cite{K87}, (2.26)-(2.28), except now we only work with components ${\cal C}$ of connected open (closed) sites that are the components of the last piece of some open (closed) path $r$ started from $\mathbf{y}$ that crosses the top, right, or bottom strips of width $n/2$ of the box $H(\mathbf{y},n)$. Denote by $a({\cal C})$ the clockwise most endpoint of such a component ${\cal C}$ such that $a({\cal C})$ still exists in $\partial ^{+}H(\mathbf{y},n)$. Given such a component ${\cal C}$, in order for it to have an $(\eta,\nu)$-fence it must satisfy (i): if ${\cal C}'$ is another such component disjoint from ${\cal C}$, then $\|a({\cal C})-a({\cal C}')\|>2\sqrt[]{\eta}2^{\nu}$, (ii): there is an appropriate crossing $s$ (vertical or horizontal as the case may be) of a small square of radius $\eta 2^{\nu}$ sitting adjacent to but outside $H(\mathbf{y},n)$ with the center of the edge adjacent to $H(\mathbf{y},n)$ equal to $a({\cal C})$, and further there is an open (closed) path from $s$ to ${\cal C}$ in a square centered at $a({\cal C})$ of the larger radius $\sqrt[]{\eta}2^{\nu}$. Here, if the component ${\cal C}$ is open and the endpoint $a({\cal C})$ is on the right side of $H(\mathbf{y},n)$ then we take $s$ to be a vertical open crossing; this crossing $s$ is the gateway to further extension of the path $r$ when the fence exists. We require further, in order that the larger squares of radii $\sqrt[]{\eta}2^{\nu}$ of condition (ii) will remain in the right half plane, that (iii): $\|a({\cal C})-(0,-n)\|>2\sqrt[]{\eta}2^{\nu}$ and $\|a({\cal C})-(0,n)\|>2\sqrt[]{\eta}2^{\nu}$. Appropriate circuits in annuli of inner and outer radii  $2\sqrt[]{\eta}2^{\nu}$ and $\sqrt[4]{\eta}2^{\nu}$ respectively and centered at $(0,\pm n)$ may be constructed with high probability to enforce condition (iii) (compare (7.5) and (7.7) in \cite{KSZ98}). Then as in \cite{K87} one may show for any $\delta >0$ there exists an $\eta>0$ such that the probability that there exists some component ${\cal C}$ that does not have an $(\eta,\nu)$-fence is at most $\delta$. This establsihes the analogue of Lemma 2 of \cite{K87} that is needed. The remainder of the proof of \eqref{E:comparePLyBPtildeLyB} follows as in the proof of Lemma 4 in \cite{K87}.  

Next it follows rather easily by Lemma 3 in \cite{K87} (extension of FKG) that \eqref{E:comparePLyBPtildeLyB} implies 
\eqref{E:comparePLyB}. Indeed, let $\mathbf{y}\in \Gamma_{H(n)}$ be given. Let $H'$ denote a vertical translation of $H(16n)$. By adjusting the vertical positioning of $H'$ appropriately, we may have for any $\mathbf{y}'\in  \Gamma_{H'}$ that $\mathbf{y}'=\mathbf{y}$. Now, independent of $\mathbf{y}'$, there is a rectangle $R'_{1}$ along the bottom of $H'$ of height at least $7n/2$ but at most $12n$ that is outside but adjacent to the bottom of $H(4n)$ and that has its top edge equal to the bottom edge of $S_{1}=S_{1}(4n)$. Trivially there is also a recatangle $R'_{2}$ of width $24n$ along the right side of $H'$ that is adjacent to but outside the right side of $H(4n)$ and that has its left side equal to the right side of $S_{2}=S_{2}(4n)$.    
Thus by the extended FKG inequality, and by using the RSW theory to estimate the probability of the existence of a vertical closed crossing of the rectangle $S_{1}\cup R'_{1}$ and also a horizontal open crossing of the rectangle $S_{2}\cup R'_{2}$ below by a constant, conditional on the event $\widetilde{{\cal L}}(\mathbf{y},H(4n))$, we may extend the closed path that $v$-tunnels through $S_{1}$ and also extend the open path that $h$-tunnels through $S_{2}$, respectively, to the bottom and right sides of $H'$ with a probability that is bounded below by a constant. Therefore by \eqref{E:comparePLyBPtildeLyB} and translation invariance we obtain \eqref{E:comparePLyB}. 

Finally we argue that \eqref{E:sizekappa2} holds. Denote $H=H(n)$ and $\widetilde{H}=H(16n)$. By \eqref{E:comparePLyB}, it follows by the disjointness of the events ${\cal L}(\widetilde{\mathbf{y}},\widetilde{H}))$, $\widetilde{\mathbf{y}}\in \Gamma_{\widetilde{H}}$, that we have
\begin{equation}\label{E:minmax}
1/(16n)\ge \min_{\widetilde{\mathbf{y}}\in \Gamma_{\widetilde{H}}}P({\cal L}(\widetilde{\mathbf{y}},\widetilde{H}))\ge (1/C)\max _{\mathbf{y}\in \Gamma_{H}}P({\cal L}(\mathbf{y},H)).
\end{equation}
But also by the RSW theory there is a constant $c>0$ such that there exists a closed horizontal crossing in each of the horizontal strips $\{n/4\le y\le n/2\}$ and $\{-n/2\le y\le -n/4\}$ in $H(n)$ as well as an open horizontal crossing of the strip $\{-n/4\le y\le n/4\}$ in $H(n)$. Therefore by disjointness of the events ${\cal L}((0,y),H(n))$, $-n/2\le y\le n/2$, we also have that
\begin{equation}\label{E:lbPLyB}
\max _{\mathbf{y}\in \Gamma_{H}}P({\cal L}(\mathbf{y},H))\ge 1/(C'n).
\end{equation}
We can now prove the statement of Lemma \ref{L:PHE3} for the special case $\rho=0$ and $\kappa=2$. Recall that $n=2^{\nu}$. We make another application of Kesten's Lemmas 2 and 4 in \cite{K87} (as slightly modified above), this time to the probability of the event ${\cal J}_{2}(0,\nu)$ such that this probability is bounded above by a constant $C_{2}$ times the probability that there exists a horizontal crossing of $H(4n)$ started from $(0,0)$. Therefore, for an upper bound we find by \eqref{E:minmax} that \[P({\cal J}_{2}(0,\nu))\le C_{2}\max_{\mathbf{y}\in \Gamma_{H(4n)}}P({\cal L}(\mathbf{y},H(4n)))\le C/n.\]  Next, by \eqref{E:minmax} and \eqref{E:lbPLyB}, by estimating a minimum from below by a constant times a maximum over a shrunken box, we have that \[P({\cal J}_{2}(0,\nu))\ge \min_{\mathbf{y}\in \Gamma_{H(n)}}P({\cal L}(\mathbf{y},H(n)))\ge 1/(C'n).\] This completes the proof of this special case. 

We now turn to the general case for $\kappa=2$. We re-orient our view to fit with the horseshoes in the statement of Lemma \ref{L:PHE3}. We assume that the origin sits in the bottom center of each of the inner and outer boxes of the pair of boxes defining $H(\rho,\nu)$, so that the horseshoe sits in the upper half-plane. For the sake of definiteness we fix disjoint intervals $\Gamma_{1}$ and $\Gamma_{2}$ of length $2^{\rho-2}$ each centered on the top and left sides of the inner horseshoe box $B_{1}(2^{\rho})$, respectively, and we take the event ${\cal J}_{2}(\rho,\nu)$ to be defined as the two-arm crossing of $H(\rho,\nu)$ such that a closed arm and an open arm meet $B_{1}(2^{\rho})$ in $\Gamma_{1}$ and $\Gamma_{2}$, respectively. Now on the event ${\cal J}_{2}(0,\nu)$ we have that ${\cal J}_{2}(0,\rho)$ occurs and also there is a two-arm path crossing the horseshoe $H(\rho,\nu)$, where this latter crossing is of course not necessarily included in the event ${\cal J}_{2}(\rho,\nu)$ as we have defined it in Section \ref{SS:multiarmcrossings}. But by Lemma 5 in \cite{K87}, the probability of this latter two arm crossing is bounded above by $CP({\cal J}_{2}(\rho,\nu))$.  Hence by the special case lower and upper bounds, we have 
\begin{equation}\label{E:lbPrhonu}
(1/C)2^{-\nu}\le P({\cal J}_{2}(0,\nu))\le CP({\cal J}_{2}(0,\rho))P({\cal J}_{2}(\rho,\nu)))\le C'2^{-\rho}P({\cal J}_{2}(\rho,\nu)).
\end{equation}
This yields the correct lower bound for $P({\cal J}_{2}(\rho,\nu))$.
To obtain an upper bound we have by Kesten's (1987) Lemma 6 that \[P({\cal J}_{2}(\rho,\nu))P({\cal J}_{2}(0,\rho))\le CP({\cal J}_{2}(0,\nu)).\] Therefore by the special case we are done with the proof of Lemma \ref{L:PHE3} for $\kappa =2$. 

We next show how to prove the Lemma \ref{L:PHE3} for $\kappa=3$. The proof follows similar lines as the proof for $\kappa=2$, but part of the trick now is to estabish an analogue of \eqref{E:lbPLyB}. This trick amounts to establishing a disjointness condition for a certain collection of half-plane three-arm events from vertices $\mathbf{y}$ (see \eqref{E:UyB} below) whose union has probability bounded below by a constant.  This idea was already explained in Zhang \cite{Z96} and mentioned besides in the Appendix A of Lawler, Schramm, and Werner \cite{LSW02b}. We assume first that $\rho=0$. As before we set $n=2^{\nu}$ and assume $\nu$ is an integer since $P({\cal J}_{3}(0,\nu))$ is decreasing in $\nu$. 
Define  
\begin{equation}\label{E:Un}
\begin{array}{ll}
{\cal U}(n):= &\exists \mbox{ a site } \mathbf{y}\in B(n/2) \mbox{ such that }\mathbf{y}\mbox{ and }\mathbf{y}+(1,0)\mbox{ are the}\\ &\mbox{unique pair of neighboring sites of highest vertical level on} \\ &\mbox{a lowest horizontal open crossing of  }B(n).
\end{array}
\end{equation}

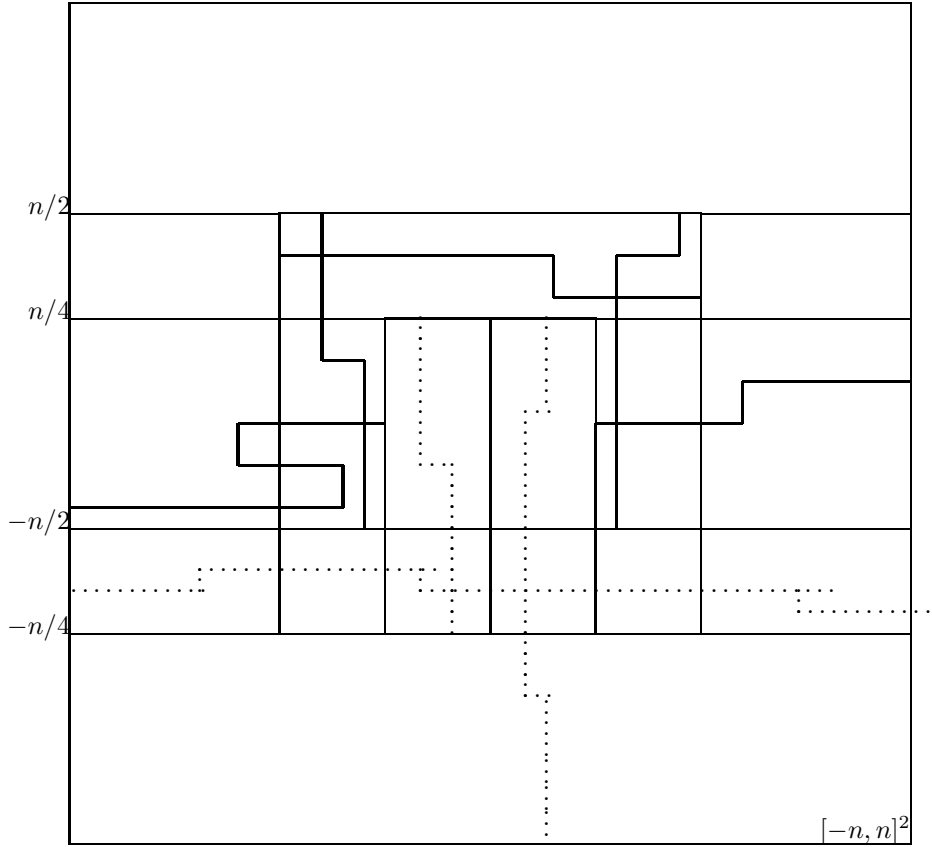
\begin{figure}\label{F:En}
\begin{center}
\setlength{\unitlength}{0.0110in}%
\begin{picture}(200,370)(67,760)

\thinlines
\put(00,800){\framebox(400,400)[br]{\mbox{$[-n, n]^2$}}}
\put(150,950){\framebox(100,100)[br]}
\put(100,900){\framebox(200,200)[br]}
\put(0,900){\line(1,0){100}}
\put(300,900){\line(1,0){100}}
\put(0,1100){\line(1,0){100}}
\put(300,1100){\line(1,0){100}}
\put(0,1050){\line(1,0){400}}
\put(0,950){\line(1,0){150}}
\put(250,950){\line(1,0){150}}
\put(150,900){\line(0,1){100}}
\put(200,900){\line(0,1){150}}
\put(250,900){\line(0,1){100}}
\put(-20, 1100){\mbox{$n/2$}}
\put(-20, 1050){\mbox{$n/4$}}
\put(-30, 900){\mbox{$-n/4$}}
\put(-30, 950){\mbox{$-n/2$}}
\put(60, 920){\mbox{${ \vdots}$}}
\put(0, 920){$\mbox{{\dots \dots \dots \dots}}$}
\put(60, 930){$\mbox{{\dots \dots \dots \dots\dots  \dots\dots}}$}
\put(165, 920){\mbox{${ \vdots}$}}
\put(165, 920){$\mbox{{\dots \dots \dots \dots\dots  \dots\dots\dots\dots\dots\dots\dots}}$}
\put(345, 910){\mbox{${ \vdots}$}}
\put(345, 910){$\mbox{{\dots \dots \dots \dots}}$}
\put(180, 900){\mbox{${ \vdots}$}}
\put(215, 890){\mbox{${ \vdots}$}}
\put(215, 880){\mbox{${ \vdots}$}}
\put(215, 870){\mbox{${ \vdots}$}}
\put(215, 870){\mbox{${ \dots}$}}
\put(225, 857){\mbox{${ \vdots}$}}
\put(225, 842){\mbox{${ \vdots}$}}
\put(225, 829){\mbox{${ \vdots}$}}
\put(225, 816){\mbox{${ \vdots}$}}
\put(225, 804){\mbox{${ \vdots}$}}
\put(180, 910){\mbox{${ \vdots}$}}
\put(180, 920){\mbox{${ \vdots}$}}
\put(180, 930){\mbox{${ \vdots}$}}
\put(180, 940){\mbox{${ \vdots}$}}
\put(180, 950){\mbox{${ \vdots}$}}
\put(180, 960){\mbox{${ \vdots}$}}
\put(180, 970){\mbox{${ \vdots}$}}
\put(165, 980){\mbox{${ \dots}$}}
\put(165, 985){\mbox{${ \vdots}$}}
\put(165, 1000){\mbox{${ \vdots}$}}
\put(165, 1015){\mbox{${ \vdots}$}}
\put(165, 1030){\mbox{${ \vdots}$}}
\put(165, 1040){\mbox{${ \vdots}$}}
\put(215, 900){\mbox{${ \vdots}$}}
\put(215, 915){\mbox{${ \vdots}$}}
\put(215, 930){\mbox{${ \vdots}$}}
\put(215, 945){\mbox{${ \vdots}$}}
\put(215, 960){\mbox{${ \vdots}$}}
\put(215, 975){\mbox{${ \vdots}$}}
\put(215, 990){\mbox{${ \vdots}$}}
\put(215, 1005){\mbox{${ \dots}$}}
\put(225, 1010){\mbox{${ \vdots}$}}
\put(225, 1025){\mbox{${ \vdots}$}}
\put(225, 1040){\mbox{${ \vdots}$}}
\thicklines
\put(0,960){\line(1,0){130}}
\put(130,960){\line(0,1){20}}
\put(130,980){\line(-1,0){50}}
\put(80,980){\line(0,1){20}}
\put(80,1000){\line(1,0){70}}
\put(250,1000){\line(1,0){70}}
\put(320,1000){\line(0,1){20}}
\put(320,1020){\line(1,0){80}}
\put(140,950){\line(0,1){80}}
\put(140,1030){\line(-1,0){20}}
\put(120,1030){\line(0,1){70}}
\put(260,950){\line(0,1){130}}
\put(260,1080){\line(1,0){30}}
\put(290,1080){\line(0,1){20}}
\put(100,1080){\line(1,0){130}}
\put(230,1080){\line(0,-1){20}}
\put(230,1060){\line(1,0){70}}
\end{picture}
\end{center}
\caption{ Paths in the event ${\cal E}(n)$. Heavy solid line: open path, dotted line: closed path. }
\end{figure}
We want to prove that $P({\cal U}(n))\ge c$. To do this we show there exists with positive probability a lowest crossing so that any vertex of a highest possible vertical level will stay in $B(n/2-1)$. This can be done as follows. See Figure 1. By independence and FKG there exists a constant $c_{1}>0$  such that $P({\cal E}(n))\ge c_{1}$ where ${\cal E}(n)$ is the joint event described by the eight conditions below: 
\begin{equation}\label{E:crossings}
\begin{array}{ll}
&\mbox{(1) }\exists \mbox{ a closed horizontal crossing of }[-n,n]\times[-n/2,0] \\ 
&\mbox{(2) }\exists \mbox{ a closed vertical crossing of }[0,n/4]\times[-n,+n/4] \\
&\mbox{(3) }\exists \mbox{ a closed vertical crossing of }[-n/4,0]\times[-n/2,+n/4] \\
&\mbox{(4) }\exists \mbox{ an open horizontal crossing of }[-n/2,n/2]\times[n/4,n/2-1] \\
&\mbox{(5) }\exists \mbox{ an open  vertical crossing of }[-n/2,-n/4]\times[0,n/2-1] \\ 
&\mbox{(6) }\exists \mbox{ an open  vertical crossing of }[n/4,n/2]\times[0,n/2-1] \\
&\mbox{(7) }\exists \mbox{ an open horizontal crossing of }[-n,-n/4]\times[0,n/4] \\
&\mbox{(8) }\exists \mbox{ an open horizontal crossing of }[n/4,n]\times[0,n/4] .
\end{array}
\end{equation}

By \eqref{E:crossings}, on ${\cal E}(n)$ there exists a lowest crossing such that a highest level on this crossing is attained for some sites $\mathbf{x}$ with $\mathbf{x}\in [-n/2,n/2]\times[n/4,n/2-1]\subset B(n/2-1)$. At the leftmost such site $\mathbf{x}$, either $\mathbf{x}$ and its immediate neighbor $\mathbf{x}+(1,0)$  already form a unique pair of highest sites on the lowest crossing or by direct surgery one may construct a unique pair of highest sites $\mathbf{y}$, $\mathbf{y}+(1,0)$ with $\mathbf{y}=\mathbf{x}+(0,1)$ by changing at most three sites in the configuration above the lowest crossing. Hence by FKG we have established that there exists a constant $c_{2}>0$ such that 
\begin{equation}\label{E:PUnEn}
P({\cal U}(n)\cap {\cal E}(n))\ge c_{2}.
\end{equation}

Now define for any $\mathbf{y}\in B(n/2)$ the event 

\begin{equation}\label{E:UyB}
\begin{array}{ll}
{\cal U}(\mathbf{y},B(n)):= &\exists \mbox{ a lowest horizontal open crossing of }B(n) \\ &\mbox{with a unique pair of highest vertices }\mathbf{y}, \mathbf{y}+(1,0).
\end{array}\end{equation}
We want to show that $P({\cal U}(\mathbf{y},B(n))\cap {\cal E}(n))$ is comparable to $1/n^{2}$ for any $\mathbf{y}\in B(n/2)$. This is done in a fashion very similar to our proof of the corresponding relation for the  case $\kappa =2$. We first establish by translation invariance and Kesten's connection method \cite{K87} that there exists a constant $C$ such that for any $\mathbf{y}\in B(n/2)$ and $\widetilde{\mathbf{y}}\in B(8n)$ we have  
\begin{equation}\label{E:comparePUyB}
P({\cal U}(\mathbf{y},B(n)))\le CP({\cal U}(\widetilde{\mathbf{y}},B(16n))).
\end{equation}
Denote by $\widetilde {{\cal U}}(\mathbf{y},B(N))$ the event that the three-arm path that issues from the highest vertex pair $\mathbf{y},\mathbf{y}+(1,0)$ in the event ${\cal U}(\mathbf{y},B(N))$ actually has the additional properties that there exist open arms that will $h$-tunnel through the squares $[-N,-3N/4]\times [-3N/4,-N/2]$ and $[3N/4,N]\times [-3N/4,-N/2]$, on the left and right sides of $B(N)$ respectively, and there exists a closed arm that will $v$-tunnel through the square $[-N/8,N/8]\times [-N,-3N/4]$ on the bottom of $B(N)$. To establish \eqref{E:comparePUyB} we shall use the proofs of Lemmas 2 and 4 in \cite{K87} to show the following.
\begin{equation}\label{E:fencearg}
P({\cal U}(\mathbf{y},B(n)))\le CP(\widetilde{{\cal U}}(\mathbf{y},B(4n))), \mbox{ uniformly for }\mathbf{y}\in B(n/2).
\end{equation}
To obtain \eqref{E:fencearg} we define the $(\eta,\nu)$-fence of \cite{K87}, Lemma 2, as above. This time we define a box $H_{3}(\mathbf{y},n)$ similar to the box $H(\mathbf{y},n)$ we worked with for the case $\kappa=2$ but oriented now with $\mathbf{y}$ in the center of the top edge of $H_{3}(\mathbf{y},n)$ instead of the center of the left edge. The formulae for the $(\eta,\nu)$-fence is re-oriented accordingly. With the new orientation we only consider components ${\cal C}$ of paths $r$ that cross strips on the right, bottom or left of $H_{3}(\mathbf{y},n)$ such that $r$ is not interrupted by the top of $H_{3}(\mathbf{y},n)$, that lies along the horizontal line $y=y_{2}$ through $\mathbf{y}$ and $\mathbf{y}+(1,0)$. If a fence exists we guarantee avoidance now of the upper left and upper right corners of $H_{3}(\mathbf{y},n)$ by endpoints $a({\cal C})$ of the corresponding components ${\cal C}$ of some arm $r$ from the pair $\mathbf{y},\mathbf{y}+(1,0)$ to the right, bottom, or left edge of $B(n)$. Define 
\begin{equation}\label{E:Kesten238}
\begin{array}{ll}
{\cal W}(\mathbf{y},B(n)):= &{\cal U}(\mathbf{y},B(n))\mbox{ occurs, and three arms that lie in } H_{3}(\mathbf{y},n)\\ &\mbox{may be chosen such that each one has an }(\eta,\nu)\mbox{-fence}.
\end{array}
\end{equation}
By construction of corridors and crossings that lie inside them that lie below the line $y=y_{2}$ and connect to the fence gates as in the proof of Lemma 4 of \cite{K87}, one obtains a constant $C(\eta)$ such that\[P({\cal W}(\mathbf{y},n))\le C(\eta)P(\widetilde{{\cal U}}(\mathbf{y},B(4n)).\] This is the appropriate analogue of inequality (2.38) in \cite{K87}. The remainder of the proof of \eqref{E:fencearg} proceeds as in \cite{K87}. 
    
To establish  \eqref{E:comparePUyB} from \eqref{E:fencearg} we proceed much as in the case $\kappa=2$. We translate the pairing of the box $B(16n)$ and the site $\widetilde{\mathbf{y}}\in B(8n)$ together to obtain a pairing of a box $B'$ of radius $16n$ and site $\mathbf{y}'$, such that $\mathbf{y}'$ coincides now with the given site $\mathbf{y}$. By \eqref{E:fencearg}, the extended FKG Lemma, RSW theory, and translation invariance we obtain  \eqref{E:comparePUyB}.   
Therefore by \eqref{E:comparePUyB} and disjointness of events we have 
\begin{equation}\label{E:Uminmax}
1/(256n^2)\ge \min_{\widetilde{\mathbf{y}}\in B(8n)}P({\cal U}(\widetilde{\mathbf{y}},B(16n)))\ge (1/C)\max _{\mathbf{y}\in B(n/2)}P({\cal U}(\mathbf{y},B(n))).
\end{equation}
Also by \eqref{E:PUnEn} and disjointness of the events ${\cal U}(\mathbf{y},B)$, $\mathbf{y}\in B(n/2)$, we have that
\begin{equation}\label{E:lbPUyB}
\max _{\mathbf{y}\in B(n/2)}P({\cal U}(\mathbf{y},B(n))\cap {\cal E}(n))\ge 1/(C'n^2).
\end{equation}
Hence by \eqref{E:Uminmax}-\eqref{E:lbPUyB} we have established that 
\begin{equation}\label{E:asympPUyB}
P({\cal U}(\mathbf{y},B(n)))\asymp 1/n^{2}, \mbox{ uniformly for } \mathbf{y}\in B(n/2).
\end{equation}

We can now prove the special case of Lemma \ref{L:PHE3} with  $\kappa=3$ when $\rho=0$. 
We make another application of the fence construction above this time to the probability of the event ${\cal J}_{3}(0,\nu)$. Therefore we obtain that this probability is bounded above by $C_{3}P({\cal U}((0,n/2),B(4n)))$. Hence, for an upper bound we find by \eqref{E:asympPUyB} that \[P({\cal J}_{3}(0,\nu))\le C_{3}\max_{\mathbf{y}\in B(2n)}P({\cal U}(\mathbf{y},B(4n)))\le C/n^2.\]  Finally, by \eqref{E:Uminmax} and \eqref{E:lbPUyB}, we have that \[P({\cal J}_{3}(0,\nu))\ge \min_{\mathbf{y}\in B(2n)}P({\cal U}(\mathbf{y},B(4n)))\ge 1/(C'n^2).\] This completes the proof of this special case. 

The proof of the general case of Lemma \ref{L:PHE3} with  $\kappa=3$ follows the same lines as the corresponding proof for $\kappa =2$. $\Box$

\section{Appendix. Proof of Lemma \ref{L:cornerprob}}\label{S:appendix}
In this section we indicate the modifications to the proof of Theorem 2 of Kesten and Zhang \cite{KZ87}, pp. 1053-1055, that are needed to establish a proof of Lemma \ref{L:cornerprob} by extending a one-arm argument to a two-arm argument.  
For the sake of concreteness we work with only two sectors, namely $S(\pi/2,n)$ and $S(\pi,n)$. These sectors are in words: the first quadrant and top half of $B(n)$, respectively. We estimate $P(\zeta(\pi/2,l,n))$ as defined in Section \ref{SS:firstmoment} directly, where for convenience we assume that $l=2^{\rho}$ and $n=2^{\nu}$ for some integers $0\le \rho<\nu$. On the event  $\zeta(\pi/2,l,n)$ there is a two-arm crossing of the annular sector $H(\pi/2,l,n):=S(\pi/2,n)\setminus S(\pi/2,l)$ (the paths end on $\partial B(n)$), where a closed path emanates from an interval $\Gamma_{1}$ of length $l/4$ along the right edge of the square sector $S(\pi/2,l)$ and an open path emanates from a corresponding fixed interval $\Gamma_{2}$ along the top edge of this square. On the event that $\zeta(\pi/2,l,n)$ occurs, we consider the clockwise-most closed path $\gamma_{1}$ as viewed from the origin that emanates from $\Gamma_{1}$ and ends on $\partial B(n)$. Since for a particular configuration of sites in the event $\gamma_{1}=r_{1}$, by changing the sub-configuration of sites sitting counterclockwise from $r_{1}$ as viewed from the origin we have another configuration of this event, we can successfully condition over all possible fixed paths  $r_{1}$ for $\gamma_{1}$. We treat $r_{1}$ as a curved boundary of each of the annular subsectors \[H(\varphi,l,n,r_{1})=\{\mathbf{x}\in H(\varphi,l,n):\mathbf{x} \mbox{ lies counterclockwise from }r_{1}\mbox{ as viewed from the origin} \}, \] $\varphi=\pi/2,\pi$. Define the half-annuli $A_{k}:=\{2^{k}\le \|(x,y)\|<2^{k+1}\mbox{ and }y\ge 0\}$, and define the events 

\begin{equation}\label{E:Fkr1}
{\cal F}_{k}(r_{1}):=\exists \mbox{ an open path in }A_{k}\cap H(\pi/2,l,n)\mbox{ from the negative }x\mbox{-axis to } r_{1}.
\end{equation}
Following \cite{KZ87} we put ${\cal E}_{j}(r_{1})={\cal F}_{3j}(r_{1})\cap {\cal F}_{3j+2}(r_{1})$, and denote by ${\cal C}_{3j}$ the innermost open path in \eqref{E:Fkr1} for $k=3j$, and by ${\cal D}_{3j+2}$ the outermostmost open path in \eqref{E:Fkr1} for $k=3j+2$. Define further $J:=\{j\ge \rho:2^{3j+2}\le n\mbox{ and }{\cal E}_{j}\mbox{ occurs}\}$ and $N:=\mbox{ cardinality of } J$. Denote by ${\cal F}(\varphi,l,n,r_{1})$ the event of an open crossing in $H(\varphi,l,n,r_{1})$ from $\Gamma_{2}$ to $\partial B(n)$. Since, given $\gamma_{1}=r_{1}$, the event $N\ge C\log (n/l)$ is increasing over the configurations of sites clockwise from $r_{1}$, we have by FKG that 
\begin{equation}\label{E:FKG}
\begin{array}{ll}
&P( {\cal F}(\pi/2,l,n,r_{1})|\gamma_{1}=r_{1})P(N\ge C\log(n/l)|\gamma_{1}=r_{1})\\ &\le P( {\cal F}(\pi/2,l,n,r_{1})\mbox{ and }N\ge C\log(n/l)|\gamma_{1}=r_{1}).
\end{array}
\end{equation}
But since, given $\gamma_{1}=r_{1}$, by RSW the event ${\cal E}_{j}$ occurs with a probability bounded below by a constant independent of $r_{1}$ and $j$, by independence we have as in \cite{KZ87} that by a simple binomial estimate $P(N\ge C\log(n/l)|\gamma_{1}=r_{1})\ge 1/2$ for a suitable constant $C$. Hence 
\begin{equation}\label{E:basicestimate}
\begin{array}{ll}
&P(\zeta(\pi/2,l,n))=\sum \limits_{r_{1}}P({\cal F}(\pi/2,l,n,r_{1})|\gamma_{1}=r_{1})P(\gamma_{1}=r_{1})
\le \\ &2\sum \limits_{r_{1}}E\left (P({\cal F}(\pi/2,l,n,r_{1})\mbox{ and }N\ge C\log (n/l)|J,{\cal C}_{3j},{\cal D}_{3j+2},j\in J)|\gamma_{1}=r_{1} \right )P(\gamma_{1}=r_{1}).
\end{array}
\end{equation}
The remainder of the proof follows just as in \cite{KZ87}, wherein one computes the conditional probability in  \eqref{E:basicestimate} over the event $N\ge C\log (n/l)$ by using the decomposition by innermost and outermost paths for a sequence $j(1)<j(2)<\cdots <j(N)$ of indices for $J$. But \[P({\cal C}_{3j}\rightarrow {\cal D}_{3j+2} \mbox{ in }H(\pi/2,l,n,r_{1})|\gamma_{1}=r_{1})\le \lambda P({\cal C}_{3j}\rightarrow {\cal D}_{3j+2} \mbox{ in }H(\pi,l,n,r_{1})|\gamma_{1}=r_{1}),\] for a suitable constant $\lambda<1$, since by a simple crossing argument and FKG, \[P({\cal C}_{3j}\rightarrow {\cal D}_{3j+2} \mbox{ in }H(\pi,l,n,r_{1})\mbox{ but not in }H(\pi/2,l,n,r_{1}))\ge c>0.\] Since also trivially, \[P({\cal D}_{3j(i)+2}\rightarrow {\cal C}_{3j(i+1)} \mbox{ in }H(\pi/2,l,n,r_{1}))\le P({\cal D}_{3j(i)+2}\rightarrow {\cal C}_{3j(i+1)} \mbox{ in }H(\pi,l,n,r_{1})),\] we obtain that 
\begin{equation}\label{E:finalestimate}
\begin{array}{ll}
P(\zeta(\pi/2,l,n))&\le 2\sum \limits_{r_{1}}\lambda^{C\log(n/l)}P({\cal F}(\pi,l,n,r_{1})\mbox{ and }N\ge C \log (n/l)|\gamma_{1}=r_{1})P(\gamma_{1}=r_{1})\\ &\le C(n/l)^{-\alpha}P(\zeta(\pi,l,n)).
\end{array}
\end{equation}
This completes the proof of Lemma \ref{L:cornerprob}.$\Box$

\section{Note Added in Proof}
The following short proof of Theorem \ref{T:minheightpivotal} based on a stationarity argument was communicated to us by Oded Schramm.

Suppose that we condition on the highest crossing $\gamma$ of the square
$[-n,n]\times[-n,n]$. If $z$ is the lowest site on $\gamma$, then an easy and
standard application of Russo-Seymour-Welsh shows that conditioned on
$\gamma$ there is probability $> \mbox{constant} > 0$ that near $z$ there is a
pivotal site (just consider closed crossings from the bottom edge of
the square to $\gamma$ near $z$).  This shows that instead of looking at
the lowest pivotal site, it is enough to consider the distance from $z$
to the bottom edge.  If $h$ is an integer in the range $[-n,0]$, then the
probability $p(h)$ that the lowest site on $\gamma$ is at height $h$ is the
same as the probability that at height $> h$ there is no crossing but at
height $\ge h$ there is a crossing. It is easy to see that $p(h)$ is
proportional to $1/n$, by the following (standard) argument.  Take
$k=n-h$.  Imagine that we take a $k \times 2n$ box (here $k$ is proportional
to $n$). If we shift down the box by $k+1$ sites, there is probability
bounded away from zero that in the new box there is a crossing while
in the old box there is none, because the boxes are disjoint. It
follows (by stationarity) that if instead we shift the box down by one
site, the probability that the new box has a crossing but the old one
does not is at least $1/(k+1)$ times a positive constant. This is
clearly a lower bound for $p(h)$.  On the other hand $p(h)$ is bounded by
the probability that there is a cluster crossing from the line $x=-n$ to
the line $x=+n$ whose lowest point on its highest crossing is at height
$h$. By stationarity, this does not depend on $h$. Since every crossing
cluster satisfies this for at most one $h$ and since the expected number
of crossing clusters meeting the horizontal strip $h-n < y < h+n$ is
bounded by a constant, we get that $p(h) \le \mbox{const}/n$ as well.

\par
\bigskip
\noindent 
{\small Dept. Mathematics, University of Colorado, Colorado Springs, CO 80933-7150.}
\par
\noindent
{\small gjmorrow@math.uccs.edu, yzhang@math.uccs.edu}
\end{document}